# The Sources of Kolmogorov's Grundbegriffe

**Glenn Shafer and Vladimir Vovk**

*Abstract.* Andrei Kolmogorov's *Grundbegriffe der Wahrscheinlichkeits-rechnung* put probability's modern mathematical formalism in place. It also provided a philosophy of probability—an explanation of how the formalism can be connected to the world of experience. In this article, we examine the sources of these two aspects of the *Grundbegriffe*—the work of the earlier scholars whose ideas Kolmogorov synthesized.

## 1. INTRODUCTION

Andrei Kolmogorov's *Grundbegriffe der Wahrscheinlichkeitsrechnung*, which set out the axiomatic basis for modern probability theory, appeared in 1933. Four years later, in his opening address to an international colloquium at the University of Geneva, Maurice Fréchet praised Kolmogorov for organizing a theory Émile Borel had created many years earlier by combining countable additivity with classical probability. Fréchet (1938b, page 54) put the matter this way in the written version of his address

> It was at the moment when Mr. Borel introduced this new kind of additivity into the calculus of probability—in 1909, that is to say—that all the elements needed to formulate explicitly the whole body of axioms of (modernized classical) probability theory came together.
>
> It is not enough to have all the ideas in mind, to recall them now and then; one must make sure that their totality is sufficient, bring them together explicitly, and take responsibility for saying that nothing further is needed in order to construct the theory.
>
> This is what Mr. Kolmogorov did. This is his achievement. (And we do not believe he wanted to claim any others, so far as the axiomatic theory is concerned.)

*Glenn Shafer is Professor, Rutgers Business School, Newark, New Jersey 07102, USA and Royal Holloway, University of London, Egham, Surrey TW20 OEX, UK (e-mail: gshafer@andromeda.rutgers.edu). Vladimir Vovk is Professor, Royal Holloway, University of London, Egham, Surrey TW20 OEX, UK (e-mail: vovk@cs.rhul.ac.uk).*

Perhaps not everyone in Fréchet's audience agreed that Borel had put everything on the table, but surely many saw the *Grundbegriffe* as a work of synthesis. In Kolmogorov's axioms and in his way of relating his axioms to the world of experience, they must have seen traces of the work of many others—the work of Borel, yes, but also the work of Fréchet himself, and that of Cantelli, Chuprov, Lévy, Steinhaus, Ulam and von Mises.

Today, what Fréchet and his contemporaries knew is no longer known. We know Kolmogorov and what came after; we have mostly forgotten what came before. This is the nature of intellectual progress, but it has left many modern students with the impression that Kolmogorov's axiomatization was born full grown—a sudden brilliant triumph over confusion and chaos.

To understand the synthesis represented by the *Grundbegriffe*, we need a broad view of the foundations of probability and the advance of measure theory from 1900 to 1930. We need to understand how measure theory became more abstract during those decades, and we need to recall what others were saying about axioms for probability, about Cournot's principle and about the relationship of probability with measure and frequency. Our review of these topics draws mainly on work by authors listed by Kolmogorov in the *Grundbegriffe*'s bibliography, especially Sergei Bernstein, Émile Borel, Francesco Cantelli, Maurice Fréchet, Paul Lévy, Antoni Łomnicki, Evgeny Slutsky, Hugo Steinhaus and Richard von Mises.

We are interested not only in Kolmogorov's mathematical formalism, but also in his philosophy of probability—how he proposed to relate the mathematical formalism to the real world. In a letter to Fréchet, Kolmogorov (1939) wrote, "You are also right in attributing to me the opinion that the formal axiomatization should be accompanied by an analysis of its real meaning." Kolmogorov devoted only two pages of the *Grundbegriffe* to such an analysis, but the question was more important to him than this brevity might suggest. We can study any mathematical formalism we like, but we have the right to call it probability only if we can explain how it relates to the phenomena classically treated by probability theory.

We begin by looking at the classical foundation that Kolmogorov's measure-theoretic foundation replaced: equally likely cases. In Section 2 we review how probability was defined in terms of equally likely cases, how the rules of the calculus of probability were derived from this definition and how this calculus was related to the real world by Cournot's principle. We also look at some paradoxes discussed at the time.

In Section 3 we sketch the development of measure theory and its increasing entanglement with probability during the first three decades of the twentieth century. This story centers on Borel, who introduced countable additivity into pure mathematics in the 1890s and then brought it to the center of probability theory, as Fréchet noted, in 1909, when he first stated and more or less proved the strong law of large numbers for coin tossing. However, the story also features Lebesgue, Radon, Fréchet, Daniell, Wiener, Steinhaus and Kolmogorov himself.

Inspired partly by Borel and partly by the challenge issued by Hilbert in 1900, a whole series of mathematicians proposed abstract frameworks for probability during the three decades we are emphasizing. In Section 4 we look at some of these, beginning with the doctoral dissertations by Rudolf Laemmel and Ugo Broggi in the first decade of the century and including an early contribution by Kolmogorov, written in 1927, five years before he started work on the *Grundbegriffe*.

In Section 5 we finally turn to the *Grundbegriffe* itself. Our review of it will confirm what Fréchet said in 1937 and what Kolmogorov says in the preface: it was a synthesis and a manual, not a report on new research. Like any textbook, its mathematics was novel for most of its readers, but its real originality was rhetorical and philosophical.

## 2. THE CLASSICAL FOUNDATION

The classical foundation of probability theory, which begins with the notion of equally likely cases, held sway for 200 years. Its elements were put in place early in the eighteenth century, and they remained in place in the early twentieth century. Even today the classical foundation is used in teaching probability.

Although twentieth century proponents of new approaches were fond of deriding the classical foundation as naive or circular, it can be defended. Its basic mathematics can be explained in a few words, and it can be related to the real world by Cournot's principle, the principle that an event with small or zero probability will not occur. This principle was advocated in France and Russia in the early years of the twentieth century, but disputed in Germany. Kolmogorov retained it in the *Grundbegriffe*.

In this section we review the mathematics of equally likely cases and recount the discussion of Cournot's principle, contrasting the support for it in France with German efforts to find other ways to relate equally likely cases to the real world. We also discuss two paradoxes, contrived at the end of the nineteenth century by Joseph Bertrand, which illustrate the care that must be taken with the concept of relative probability. The lack of consensus on how to make philosophical sense of equally likely cases and the confusion revealed by Bertrand's paradoxes were two sources of dissatisfaction with the classical theory.

### 2.1 The Classical Calculus

The classical definition of probability was formulated by Jacob Bernoulli (1713) in *Ars Conjectandi* and Abraham de Moivre in 1718 in *The Doctrine of Chances*: the probability of an event is the ratio of the number of equally likely cases that favor it to the total number of equally likely cases possible under the circumstances.

From this definition, de Moivre derived two rules for probability. The *theorem of total probability*, or the *addition theorem*, says that if $A$ and $B$ cannot both happen, then

$$\text{probability of } A \text{ or } B \text{ happening}$$

$$= \frac{\# \text{ of cases favoring } A \text{ or } B}{\text{total } \# \text{ of cases}}$$

$$= \frac{\# \text{ of cases favoring } A}{\text{total } \# \text{ of cases}} + \frac{\# \text{ of cases favoring } B}{\text{total } \# \text{ of cases}}$$

$$= (\text{probability of } A) + (\text{probability of } B).$$

The *theorem of compound probability*, or the *multiplication theorem*, says

$$\text{probability of both } A \text{ and } B \text{ happening}$$

$$= \frac{\# \text{ of cases favoring both } A \text{ and } B}{\text{total } \# \text{ of cases}}$$

$$= \frac{\# \text{ of cases favoring } A}{\text{total } \# \text{ of cases}}$$

$$\cdot \frac{\text{\# of cases favoring both } A \text{ and } B}{\text{\# of cases favoring } A}$$

$$= (\text{probability of } A)$$

$$\cdot (\text{probability of } B \text{ if } A \text{ happens}).$$

These arguments were still standard fare in probability textbooks at the beginning of the twentieth century, including the great treatises by Henri Poincaré (1896) in France, Andrei Markov (1900) in Russia and Emanuel Czuber (1903) in Germany. Some years later we find them in Guido Castelnuovo's (1919) Italian textbook, which has been held out as the acme of the genre (Onicescu, 1967).

Geometric probability was incorporated into the classical theory in the early nineteenth century. Instead of counting equally likely cases, one measures their geometric extension—their area or volume. However, probability is still a ratio, and the rules of total and compound probability are still theorems. This was explained by Antoine-Augustin Cournot (1843, page 29) in his influential treatise on probability and statistics, *Exposition de la théorie des chances et des probabilités*. This understanding of geometric probability did not change in the early twentieth century, when Borel and Lebesgue expanded the class of sets for which we can define geometric extension. We may now have more events with which to work, but we define and study geometric probabilities as before. Cournot would have seen nothing novel in Felix Hausdorff's (1914, pages 416–417) definition of probability in the chapter on measure theory in his treatise on set theory.

The classical calculus was enriched at the beginning of the twentieth century by a formal and universal notation for relative probabilities. Hausdorff (1901) introduced the symbol $p_F(E)$ for what he called the *relative Wahrscheinlichkeit von $E$, posito $F$* (relative probability of $E$ given $F$). Hausdorff explained that this notation can be used for any two events $E$ and $F$, no matter what their temporal or logical relationship, and that it allows one to streamline Poincaré's proofs of the addition and multiplication theorems. Hausdorff's notation was adopted by Czuber in 1903. Kolmogorov used it in the *Grundbegriffe*, and it persisted, especially in the German literature, until the middle of the twentieth century, when it was displaced by the more flexible $P(E|F)$, which Harold Jeffreys (1931) introduced in his *Scientific Inference*.

## 2.2 Cournot's Principle

An event with very small probability is *morally impossible*: it will not happen. Equivalently, an event with very high probability is *morally certain*: it will happen. This principle was first formulated within mathematical probability by Jacob Bernoulli. In his *Ars Conjectandi*, published in 1713, Bernoulli proved a celebrated theorem: in a sufficiently long sequence of independent trials of an event, there is a very high probability that the frequency with which the event happens will be close to its probability. Bernoulli explained that we can treat the very high probability as moral certainty and so use the frequency of the event as an estimate of its probability.

Probabilistic moral certainty was widely discussed in the eighteenth century. In the 1760s, the French savant Jean d'Alembert muddled matters by questioning whether the prototypical event of very small probability, a long run of many happenings of an event as likely to fail as happen on each trial, is possible at all. A run of a hundred may be metaphysically possible, he felt, but it is physically impossible. It has never happened and never will happen (d'Alembert, 1761, 1767; Daston, 1979). Buffon (1777) argued that the distinction between moral and physical certainty is one of degree. An event with probability 9999/10,000 is morally certain; an event with much greater probability, such as the rising of the sun, is physically certain (Loveland, 2001).

Cournot, a mathematician now remembered as an economist and a philosopher of science (Martin, 1996, 1998), gave the discussion a nineteenth century cast. Being equipped with the idea of geometric probability, Cournot could talk about probabilities that are vanishingly small. He brought physics to the foreground. It may be mathematically possible, he argued, for a heavy cone to stand in equilibrium on its vertex, but it is physically impossible. The event's probability is vanishingly small. Similarly, it is physically impossible for the frequency of an event in a long sequence of trials to differ substantially from the event's probability (Cournot, 1843, pages 57 and 106).

In the second half of the nineteenth century, the principle that an event with a vanishingly small probability will not happen took on a real role in physics, most saliently in Ludwig Boltzmann's statistical understanding of the second law of thermodynamics. As Boltzmann explained in the 1870s, dissipative processes are irreversible because the probability of a state with entropy far from the maximum is vanishingly small (von Plato, 1994, page 80; Seneta, 1997). Also notable was Henri Poincaré's use of the principle in celestial mechanics. Poincaré's (1890) recurrence theorem says that an isolated mechanical system confined to a bounded region of its phase space will eventually return arbitrarily close to its initial state, provided only that this initial state is not exceptional. The states for which the recurrence does not hold are exceptional inasmuch as they are contained in subregions whose total volume is arbitrarily small.

Saying that an event of very small or vanishingly small probability will not happen is one thing. Saying that probability theory gains empirical meaning only by ruling out the happening of such events is another. Cournot (1843, page 78) seems to have been the first to say explicitly that probability theory does gain empirical meaning only by declaring events of vanishingly small probability to be impossible:

> ... *The physically impossible event is therefore the one that has infinitely small probability*, and only this remark gives substance—objective and phenomenal value—to the theory of mathematical probability.

[The phrase "objective and phenomenal" refers to Kant's distinction between the noumenon, or thing-in-itself, and the phenomenon, or object of experience (Daston, 1994).] After the Second World War, some authors began to use "Cournot's principle" for the principle that an event of very small or zero probability singled out in advance will not happen, especially when this principle is advanced as the unique means by which a probability model is given empirical meaning.

2.2.1 *The viewpoint of the French probabilists.* In the early decades of the twentieth century, probability theory was beginning to be understood as pure mathematics. What does this pure mathematics have to do with the real world? The mathematicians who revived research in probability theory in France during these decades, Émile Borel, Jacques Hadamard, Maurice Fréchet and Paul Lévy, made the connection by treating events of small or zero probability as impossible.

Borel explained this repeatedly, often in a style more literary than mathematical or philosophical (Borel, 1906, 1909b, 1914, 1930). Borel's many discussions of the considerations that go into assessing the boundaries of practical certainty culminated in a classification more refined than Buffon's. A probability of $10^{-6}$, he decided, is negligible at the human scale, a probability of $10^{-15}$ at the terrestrial scale and a probability of $10^{-50}$ at the cosmic scale (Borel, 1939, pages 6–7).

Hadamard, the preeminent analyst who did pathbreaking work on Markov chains in the 1920s (Bru, 2003), made the point in a different way. Probability theory, he said, is based on two notions: the notion of perfectly equivalent (equally likely) events and the notion of a very unlikely event (Hadamard, 1922, page 289). Perfect equivalence is a mathematical assumption which cannot be verified. In practice, equivalence is not perfect—one of the grains in a cup of sand may be more likely than another to hit the ground first when they are thrown out of the cup—but this need not prevent us from applying the principle of the very unlikely event. Even if the grains are not exactly the same, the probability of any particular one hitting the ground first is negligibly small. Hadamard was the teacher of both Fréchet and Lévy.

Among the French mathematicians of this period, it was Lévy who expressed most clearly the thesis that Cournot's principle is probability's only bridge to reality. In his *Calcul des probabilités* (Lévy, 1925) Lévy emphasized the different roles of Hadamard's two notions. The notion of equally likely events, Lévy explained, suffices as a foundation for the mathematics of probability, but so long as we base our reasoning only on this notion, our probabilities are merely subjective. It is the notion of a very unlikely event that permits the results of the mathematical theory to take on practical significance (Lévy, 1925, pages 21, 34; see also Lévy, 1937, page 3). Combining the notion of a very unlikely event with Bernoulli's theorem, we obtain the notion of the objective probability of an event, a physical constant that is measured by frequency. Objective probability, in Lévy's view, is entirely analogous to length and weight, other physical constants whose empirical meaning is also defined by methods established for measuring them to a reasonable approximation (Lévy, 1925, pages 29–30).

By the time he undertook to write the *Grundbegriffe*, Kolmogorov must have been very familiar with Lévy's views. He had cited Lévy's 1925 book in his 1931 article on Markov processes and subsequently, during his visit to France, had spent a great deal of time talking with Lévy about probability. He could also have learned about Cournot's principle from the Russian literature. The champion of the principle in Russia had been Chuprov, who became professor of statistics in Petersburg in 1910. Chuprov put Cournot's principle—which he called Cournot's lemma—at the heart of this project; it was, he said, a basic principle of the logic of the probable (Chuprov, 1910; Sheynin, 1996, pages 95–96). Markov, who also worked in Petersburg, learned about the burgeoning field of mathematical statistics from Chuprov (Ondar, 1981), and we see an echo of Cournot's principle in Markov's (1912, page 12 of the German edition) textbook:

> The closer the probability of an event is to one, the more reason we have to expect the event to happen and not to expect its opposite to happen.
>
> In practical questions, we are forced to regard as certain events whose probability comes more or less close to one, and to regard as impossible events whose probability is small.
>
> Consequently, one of the most important tasks of probability theory is to identify those events whose probabilities come close to one or zero.

The Russian statistician Evgeny Slutsky discussed Chuprov's views in his influential article on limit theorems (Slutsky, 1925). Kolmogorov included Lévy's book and Slutsky's article in his bibliography, but not Chuprov's book. An opponent of the Bolsheviks, Chuprov was abroad when they seized power, and he never returned home. He remained active in Sweden and Germany, but his health soon failed, and he died in 1926 at the age of 52.

2.2.2 *Strong and weak forms of Cournot's principle.* Cournot's principle has many variations. Like probability, moral certainty can be subjective or objective. Some authors make moral certainty sound truly equivalent to absolute certainty; others emphasize its pragmatic meaning.

For our story, it is important to distinguish between the strong and weak forms of the principle (Fréchet, 1951, page 6; Martin, 2003). The strong form refers to an event of small or zero probability that we single out in advance of a single trial: it says the event will not happen on that trial. The weak form says that an event with very small probability will happen very rarely in repeated trials.

Borel, Lévy and Kolmogorov all subscribed to Cournot's principle in its strong form. In this form, the principle combines with Bernoulli's theorem to produce the unequivocal conclusion that an event's probability will be approximated by its frequency in a particular sufficiently long sequence of independent trials. It also provides a direct foundation for statistical testing. If the meaning of probability resides precisely in the nonhappening of small-probability events singled out in advance, then we need no additional principles to justify rejecting a hypothesis that gives small probability to an event we single out in advance and then observe to happen.

Other authors, including Chuprov, enunciated Cournot's principle in its weak form, and this can lead in a different direction. The weak principle combines with Bernoulli's theorem to produce the conclusion that an event's probability will usually be approximated by its frequency in a sufficiently long sequence of independent trials, a general principle that has the weak principle as a special case. This was pointed out in the famous textbook by Castelnuovo (1919, page 108). On page 3, Castelnuovo called the general principle the *empirical law of chance*:

> In a series of trials repeated a large number of times under identical conditions, each of the possible events happens with a (relative) frequency that gradually equals its probability. The approximation usually improves as the number of trials increases.

Although the special case where the probability is close to 1 is sufficient to imply the general principle, Castelnuovo preferred to begin his introduction to the meaning of probability by enunciating the general principle, and so he can be considered a frequentist. His approach was influential. Maurice Fréchet and Maurice Halbwachs adopted it in their textbook in 1924. It brought Fréchet to the same understanding of objective probability as Lévy: objective probability is a physical constant that is measured by frequency (Fréchet, 1938a, page 5; 1938b, pages 45–46).

The weak point of Castelnuovo and Fréchet's position lies in the modesty of their conclusion: they conclude only that an event's probability is *usually* approximated by its frequency. When we estimate a probability from an observed frequency, we are taking a further step: we are assuming that what usually happens has happened in the particular case. This step requires the strong form of Cournot's principle. According to Kolmogorov (1956, page 240 of the 1965 English edition), it is a reasonable step only if we have some reason to assume that the position of the particular case among other potential ones "is a regular one, that is, that it has no special features."

2.2.3 *British indifference and German skepticism.* The mathematicians who worked on probability in France in the early twentieth century were unusual in the extent to which they delved into the philosophical side of their subject. Poincaré had made a mark in the philosophy of science as well as in mathematics, and Borel, Fréchet and Lévy tried to emulate him. The situation in Britain and Germany was different.

In Britain there was little mathematical work in probability proper in this period. In the nineteenth century, British interest in probability had been practical and philosophical, not mathematical (Porter, [1986](#), page 74ff). Robert Leslie Ellis ([1849](#)) and John Venn ([1888](#)) accepted the usefulness of probability, but insisted on defining it directly in terms of frequency, leaving no role for Bernoulli's theorem and Cournot's principle (Daston, [1994](#)). These attitudes persisted even after Pearson and Fisher brought Britain into a leadership role in mathematical statistics. The British statisticians had no puzzle to solve concerning how to link probability to the world. They were interested in reasoning directly about frequencies.

In contrast with Britain, Germany did see a substantial amount of mathematical work in probability during the first decades of the twentieth century, much of it published in German by Scandinavians and eastern Europeans, but few German mathematicians of the first rank fancied themselves philosophers. The Germans were already pioneering the division of labor to which we are now accustomed, between mathematicians who prove theorems about probability, and philosophers, logicians, statisticians and scientists who analyze the meaning of probability. Many German statisticians believed that one must decide what level of probability will count as practical certainty in order to apply probability theory (von Bortkiewicz, [1901](#), page 825; Bohlmann, [1901](#), page 861), but German philosophers did not give Cournot's principle a central role.

The most cogent and influential of the German philosophers who discussed probability in the late nineteenth century was Johannes von Kries ([1886](#)), whose *Principien der Wahrscheinlichkeitsrechnung* first appeared in 1886. von Kries rejected what he called the orthodox philosophy of Laplace and the mathematicians who followed him. As von Kries saw it, these mathematicians began with a subjective concept of probability, but then claimed to establish the existence of objective probabilities by means of a so-called law of large numbers, which they erroneously derived by combining Bernoulli's theorem with the belief that small probabilities can be neglected. Having both subjective and objective probabilities at their disposal, these mathematicians then used Bayes' theorem to reason about objective probabilities for almost any question where many observations are available. All this, von Kries believed, was nonsense. The notion that an event with very small probability is impossible was, in von Kries' eyes, simply d'Alembert's mistake.

von Kries believed that objective probabilities sometimes exist, but only under conditions where equally likely cases can legitimately be identified. Two conditions, he thought, are needed:

- Each case is produced by equally many of the possible arrangements of the circumstances, and this remains true when we look back in time to earlier circumstances that led to the current ones. In this sense, the relative sizes of the cases are *natural*.
- Nothing besides these circumstances affects our expectation about the cases. In this sense, the *Spielräume* are *insensitive.* [In German, *Spiel* means game or play, and *Raum* (plural Räume) means room or space. In most contexts, *Spielraum* can be translated as leeway or room for maneuver. For von Kries the *Spielraum* for each case was the set of all arrangements of the circumstances that produce it.]

von Kries' *principle of the Spielräume* was that objective probabilities can be calculated from equally likely cases when these conditions are satisfied. He considered this principle analogous to Kant's principle that everything that exists has a cause. Kant thought that we cannot reason at all without the principle of cause and effect. von Kries thought that we cannot reason about objective probabilities without the principle of the *Spielräume*.

Even when an event has an objective probability, von Kries saw no legitimacy in the law of large numbers. Bernoulli's theorem is valid, he thought, but it tells us only that a large deviation of an event's frequency from its probability is just as unlikely as some other unlikely event, say a long run of successes. What will actually happen is another matter. This disagreement between Cournot and von Kries can be seen as a quibble about words. Do we say that an event will not happen (Cournot) or do we say merely that it is as unlikely as some other event we do not expect to happen (von Kries)? Either way, we proceed as if it will not happen. However, the quibbling has its reasons. Cournot wanted to make a definite prediction, because this provides a bridge from probability theory to the world of phenomena—the real world, as those who have not studied Kant would say. von Kries thought he had a different way to connect probability theory with phenomena.

von Kries' critique of moral certainty and the law of large numbers was widely accepted in Germany (Kamlah, 1983). Czuber, in the influential textbook we have already mentioned, named Bernoulli, d'Alembert, Buffon and De Morgan as advocates of moral certainty and declared them all wrong; the concept of moral certainty, he said, violates the fundamental insight that an event of ever so small a probability can still happen (Czuber, 1843, page 15; see also Meinong, 1915, page 591).

This wariness about ruling out the happening of events whose probability is merely very small does not seem to have prevented acceptance of the idea that zero probability represents impossibility. Beginning with Wiman's work on continued fractions in 1900, mathematicians writing in German worked on showing that various sets have measure zero, and everyone understood that the point was to show that these sets are impossible (see Felix Bernstein, 1912, page 419). This suggests a great gulf between zero probability and merely small probability. One does not sense such a gulf in the writings of Borel and his French colleagues; as we have seen, the vanishingly small, for them, was merely an idealization of the very small.

von Kries' principle of the *Spielräume* did not endure, because no one knew how to use it, but his project of providing a Kantian justification for the uniform distribution of probabilities remained alive in German philosophy in the first decades of the twentieth century (Meinong, 1915; Reichenbach, 1916). John Maynard Keynes (1921) brought it into the English literature, where it continues to echo, to the extent that today's probabilists, when asked about the philosophical grounding of the classical theory of probability, are more likely to think about arguments for a uniform distribution of probabilities than about Cournot's principle.

### 2.3 Bertrand's Paradoxes

How do we know cases are equally likely, and when something happens, do the cases that remain possible remain equally likely? In the decades before the *Grundbegriffe*, these questions were frequently discussed in the context of paradoxes formulated by Joseph Bertrand, an influential French mathematician, in a textbook published in 1889.

We now look at discussions by other authors of two of Bertrand's paradoxes: Poincaré's discussion of the paradox of the three jewelry boxes and Borel's discussion of the paradox of the great circle. (In the literature of the period, "Bertrand's paradox" usually referred to a third paradox, concerning two possible interpretations of the idea of choosing a random chord on a circle. Determining a chord by choosing two random points on the circumference is not the same as determining it by choosing a random distance from the center and then a random orientation.) The paradox of the great circle was also discussed by Kolmogorov and is now sometimes called the *Borel–Kolmogorov paradox*.

2.3.1 *The paradox of the three jewelry boxes.* This paradox, laid out by Bertrand (1889, pages 2–3), involves three identical jewelry boxes, each with two drawers. Box A has gold medals in both drawers, box B has silver medals in both, and box C has a gold medal in one and a silver medal in the other. Suppose we choose a box at random. It will be box C with probability 1/3. Now suppose we open at random one of the drawers in the box we have chosen. There are two possibilities for what we find:

- We find a gold medal. In this case, only two possibilities remain: the other drawer has a gold medal (we have chosen box A) or the other drawer has a silver medal (we have chosen box C).
- We find a silver medal. Here also, only two possibilities remain: the other drawer has a gold medal (we have chosen box C) or the other drawer has a silver medal (we have chosen box B).

Either way, it seems, there are now two cases, one of which is that we have chosen box C. So the probability that we have chosen box C is now 1/2.

Bertrand himself did not accept the conclusion that opening the drawer would change the probability of having box C from 1/3 to 1/2, and Poincaré (1912, pages 26–27) gave an explanation: Suppose the drawers in each box are labeled (where we cannot see) $\alpha$ and $\beta$, and suppose the gold medal in box C is in drawer $\alpha$. Then there are six equally likely cases for the drawer we open:

1. Box A, drawer $\alpha$: gold medal.
2. Box A, drawer $\beta$: gold medal.
3. Box B, drawer $\alpha$: silver medal.
4. Box B, drawer $\beta$: silver medal.
5. Box C, drawer $\alpha$: gold medal.
6. Box C, drawer $\beta$: silver medal.

When we find a gold medal, say, in the drawer we have opened, three of these cases remain possible: case 1, case 2 and case 5. Of the three, only one favors our having our hands on box C, so the probability for box C is still 1/3.

2.3.2 *The paradox of the great circle.* Bertrand (1889, pages 6–7) begins with a simple question: if we choose at random two points on the surface of a sphere, what is the probability that the distance between them is less than $10'$?

By symmetry, we can suppose that the first point is known. So one way to answer the question is to calculate the proportion of a sphere's surface that lies within $10'$ of a given point. This is $2.1 \times 10^{-6}$.

Bertrand also found a different answer. After fixing the first point, he said, we can also assume that we know the great circle that connects the two points, because the possible chances are the same on great circles through the first point. There are 360 degrees—2160 arcs of $10'$ each—in this great circle. Only the points in the two neighboring arcs are within $10'$ of the first point, and so the probability sought is $2/2160$, or $9.3 \times 10^{-4}$. This is many times larger than the probability found by the first method. Bertrand considered both answers equally valid, the original question being ill-posed. The concept of choosing points at random on a sphere was not, he said, sufficiently precise.

In his own probability textbook Borel (1909b, pages 100–104) explained that Bertrand was mistaken. Bertrand's first method, based on the assumption that equal areas on the sphere have equal chances of containing the second point, is correct. His second method, based on the assumption that equal arcs on a great circle have equal chances of containing it, is incorrect. Writing M and M′ for the two points to be chosen at random on the sphere, Borel explained Bertrand's mistake as follows:

> . . . The error begins when, after fixing the point M and the great circle, one assumes that the probability of M′ being on a given arc of the great circle is proportional to the length of that arc. If the arcs have no width, then in order to speak rigorously, we must assign the value zero to the probability that M and M′ are on the circle. In order to avoid this factor of zero, which makes any calculation impossible, one must consider a thin bundle of great circles all going through M, and then it is obvious that there is a greater probability for M′ to be situated in a vicinity 90 degrees from M than in the vicinity of M itself (Fig. 13).

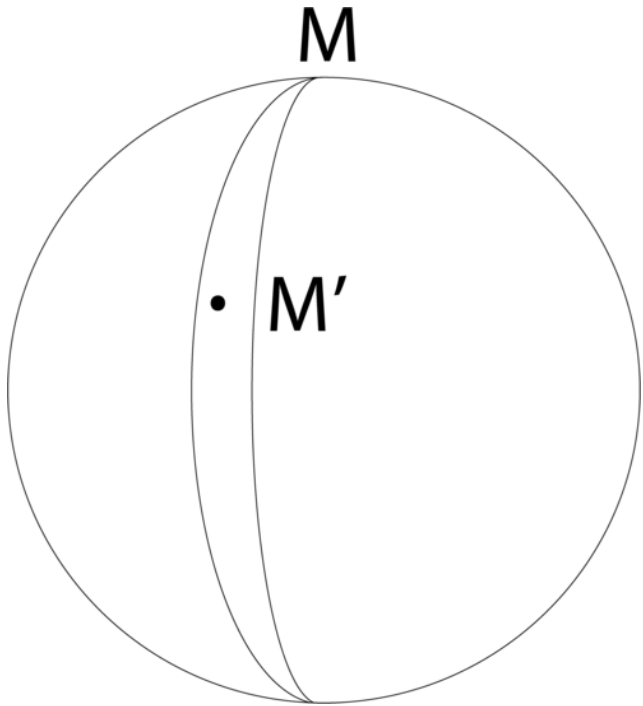

Fig. 1. *Borel's Figure* 13.

To give this argument practical content, Borel discussed how one might measure the longitude of a point on the surface of the earth. If we use astronomical observations, then we are measuring an angle, and errors in the measurement of the angle correspond to wider distances on the ground at the equator than at the poles. If we instead use geodesic measurements, say with a line of markers on each of many meridians, then to keep the markers out of each other's way, we must make them thinner and thinner as we approach the poles.

2.3.3 *Appraisal.* Poincaré, Borel and others who understood the principles of the classical theory were able to resolve the paradoxes that Bertrand contrived. Two principles emerge from the resolutions they offered:

- The equally likely cases must be detailed enough to represent new information (e.g., we find a gold medal) in all relevant detail. The remaining equally likely cases will then remain equally likely.
- We may need to consider the real observed event of nonzero probability that is represented in an idealized way by an event of zero probability (e.g., a randomly chosen point falls on a particular meridian). We should pass to the limit only after absorbing the new information.

Not everyone found it easy to apply these principles, however, and the confusion surrounding the paradoxes was another source of dissatisfaction with the classical theory.

## 3. MEASURE-THEORETIC PROBABILITY BEFORE THE GRUNDBEGRIFFE

A discussion of the relationship between measure and probability in the first decades of the twentieth century must navigate many pitfalls, because measure theory itself evolved, beginning as a theory about the measurability of sets of real numbers and then becoming more general and abstract. Probability theory followed along, but since the meaning of *measure* was changing, we can easily misunderstand things said at the time about the relationship between the two theories.

The development of theories of measure and integration during the late nineteenth and early twentieth centuries has been studied extensively (Hawkins, 1975; Pier, 1994a). Here we offer only a bare-bones sketch, beginning with Borel and Lebesgue, and touching on those steps that proved most significant for the foundations of probability. We discuss the work of Carathéodory, Radon, Fréchet and Nikodym, who made measure primary and integral secondary, as well as the contrasting approach of Daniell, who took integration to be basic, and Wiener, who applied Daniell's methods to Brownian motion. Then we discuss Borel's strong law of large numbers, which focused attention on measure rather than on integration. After looking at Steinhaus' axiomatization of Borel's denumerable probability, we turn to Kolmogorov's use of measure theory in probability in the 1920s.

### 3.1 Measure Theory from Borel to Fréchet

Émile Borel is considered the founder of measure theory. Whereas Peano and Jordan had extended the concept of length from intervals to a larger class of sets of real numbers by approximating the sets inside and outside with finite unions of intervals, Borel used countable unions. His motivation came from complex analysis. In his doctoral dissertation Borel (1895) studied certain series that were known to diverge on a dense set of points on a closed curve and hence, it was thought, could not be continued analytically into the region bounded by the curve. Roughly speaking, Borel discovered that the set of points where divergence occurred, although dense, can be covered by a countable number of intervals with arbitrarily small total length. Elsewhere on the curve—almost everywhere, we would say now—the series does converge and so analytic continuation is possible (Hawkins, 1975, Section 4.2). This discovery led Borel to a new theory of measurability for subsets of $[0, 1]$ (Borel, 1898).

Borel's innovation was quickly seized upon by Henri Lebesgue, who made it the basis for his powerful theory of integration (Lebesgue, 1901). We now speak of Lebesgue measure on the real numbers $R$ and on the $n$-dimensional space $R^n$, and of the Lebesgue integral in these spaces. We need not review Lebesgue's theory, but we should mention one theorem, the precursor of the Radon–Nikodym theorem: any countably additive and absolutely continuous set function on the real numbers is an indefinite integral. This result first appeared in (Lebesgue, 1904; Hawkins, 1975, page 145; Pier, 1994a, page 524). He generalized it to $R^n$ in 1910 (Hawkins, 1975, page 186).

Wacław Sierpiński (1918) gave an axiomatic treatment of Lebesgue measure. In this note, important to us because of the use Hugo Steinhaus later made of it, Sierpiński characterized the class of Lebesgue measurable sets as the smallest class $K$ of sets that satisfy the following conditions:

I. For every set $E$ in $K$, there is a nonnegative number $\mu(E)$ that will be its measure and will satisfy conditions II, III, IV and V.
II. Every finite closed interval is in $K$ and has its length as its measure.
III. The class $K$ is closed under finite and countable unions of disjoint elements, and $\mu$ is finitely and countably additive.
IV. If $E_1 \supset E_2$, and $E_1$ and $E_2$ are in $K$, then $E_1 \setminus E_2$ is in $K$.
V. If $E$ is in $K$ and $\mu(E) = 0$, then any subset of $E$ is in $K$.

An arbitrary class $K$ that satisfies these conditions is not necessarily a field; there is no requirement that the intersection of two of $K$'s elements also be in $K$.

Lebesgue's measure theory was first made abstract by Johann Radon (1913). Radon unified Lebesgue and Stieltjes integration by generalizing integration with respect to Lebesgue measure to integration with respect to any countably additive set function on the Borel sets in $R^n$. The generalization included a version of the theorem of Lebesgue we just mentioned: if a countably additive set function $g$ on $R^n$ is absolutely continuous with respect to another countably additive set function $f$, then $g$ is an indefinite integral with respect to $f$ (Hawkins, 1975, page 189).

Constantin Carathéodory was also influential in drawing attention to measures on Euclidean spaces other than Lebesgue measure. Carathéodory (1914) gave axioms for outer measure in a $q$-dimensional space, derived the notion of measure and applied these ideas not only to Lebesgue measure on Euclidean spaces, but also to lower dimensional measures on Euclidean space which assign lengths to curves, areas to surfaces and so forth (Hochkirchen, 1999). Carathéodory also recast Lebesgue's theory of integration to make measure even more fundamental; in his textbook (Carathéodory, 1918) on real functions, he defined the integral of a positive function on a subset of $R^n$ as the $(n + 1)$-dimensional measure of the region between the subset and the function's graph (Bourbaki, 1994, page 228).

It was Fréchet who first went beyond Euclidean space. Fréchet (1915a, b) observed that much of Radon's reasoning does not depend on the assumption that one is working in $R^n$. One can reason in the same way in a much larger space, such as a space of functions. Any space will do, so long as the countably additive set function is defined on a $\sigma$-field of its subsets, as Radon had required. Fréchet did not, however, manage to generalize Radon's theorem on absolute continuity to the fully abstract framework. This generalization, now called the *Radon–Nikodym theorem*, was obtained by Otton Nikodym fifteen years later (Nikodym, 1930).

Did Fréchet himself have probability in mind when he proposed a calculus that allows integration over function space? Probably so. An integral is a mean value. In a Euclidean space this might be a mean value with respect to a distribution of mass or electrical charge, but we cannot distribute mass or charge over a space of functions. The only thing we can imagine distributing over such a space is probability or frequency. However, Fréchet thought of probability as an application of mathematics, not as a branch of pure mathematics itself, so he did not think he was axiomatizing probability theory.

It was Kolmogorov who first called Fréchet's theory a foundation for probability theory. He put the matter this way in the preface to the *Grundbegriffe*:

> . . . After Lebesgue's investigations, the analogy between the measure of a set and the probability of an event, as well as between the integral of a function and the mathematical expectation of a random variable, was clear. This analogy could be extended further; for example, many properties of independent random variables are completely analogous to corresponding properties of orthogonal functions. But in order to base probability theory on this analogy, one still needed to liberate the theory of measure and integration from the geometric elements still in the foreground with Lebesgue. This liberation was accomplished by Fréchet.

It should not be inferred from this passage that Fréchet and Kolmogorov used "measure" in the way we do today. Fréchet may have liberated measure and integration from its geometric roots, but Fréchet and Kolmogorov continued to reserve the word measure for geometric settings. Throughout the 1930s, what we now call a measure, they called an *additive set function*. The usage to which we are now accustomed became standard only after the Second World War.

### 3.2 Daniell's Integral and Wiener's Differential Space

Percy Daniell, an Englishman working at the Rice Institute in Houston, Texas, introduced his integral in a series of articles (Daniell, 1918, 1919a, b, 1920) in the *Annals of Mathematics*.

Like Fréchet, Daniell considered an abstract set $E$, but instead of beginning with an additive set function on subsets of $E$, he began with what he called an integral on $E$—a linear operator on some class $T_0$ of real-valued functions on $E$. The class $T_0$ might consist of all continuous functions (if $E$ is endowed with a topology) or perhaps all step functions. Applying Lebesgue's methods in this general setting, Daniell extended the linear operator to a wider class $T_1$ of functions on $E$, the summable functions. In this way, the Riemann integral is extended to the Lebesgue integral, the Stieltjes integral is extended to the Radon integral and so on (Daniell, 1918). Using ideas from Fréchet's dissertation, Daniell also gave examples in infinite-dimensional spaces (Daniell, 1919a, b). Daniell (1921) even used his theory of integration to construct a theory of Brownian motion. However, he did not succeed in gaining recognition for this last contribution; it seems to have been completely ignored until Stephen Stigler spotted it in the 1970s (Stigler, 1973).

The American ex-child prodigy and polymath Norbert Wiener, when he came upon Daniell's 1918 and July 1919 articles (Daniell, 1918, 1919a), was in a better position than Daniell himself to appreciate and advertise their remarkable potential for probability (Wiener, 1956; Masani, 1990). Having studied philosophy as well as mathematics, Wiener was well aware of the intellectual significance of Brownian motion and of Einstein's mathematical model for it.

In November 1919, Wiener submitted his first article (Wiener, 1920) on Daniell's integral to the *Annals of Mathematics*, the journal where Daniell's four articles on it had appeared. This article did not yet discuss Brownian motion; it merely laid out a general method for setting up a Daniell integral when the underlying space $E$ is a function space. However, by August 1920, Wiener was in France to explain his ideas on Brownian motion to Fréchet and Lévy (Segal, 1992, page 397). He followed up with a series of articles (Wiener, 1921a, b), including a later much celebrated article on "differential-space" (Wiener, 1923).

Wiener's basic idea was simple. Suppose we want to formalize the notion of Brownian motion for a finite time interval, say $0 \leq t \leq 1$. A realized path is a function on $[0, 1]$. We want to define mean values for certain functionals (real-valued functions of the realized path). To set up a Daniell integral that gives these mean values, Wiener took $T_0$ to consist of functionals that depend only on the path's values at a finite number of time points. One can find the mean value of such a functional using Gaussian probabilities for the changes from each time point to the next. Extending this integral by Daniell's method, he succeeded in defining mean values for a wide class of functionals. In particular, he obtained probabilities (mean values for indicator functions) for certain sets of paths. He showed that the set of continuous paths has probability 1, while the set of differentiable paths has probability 0.

It is now commonplace to translate this work into Kolmogorov's measure-theoretic framework. Kiyoshi Itô, for example, in a commentary published along with Wiener's articles from this period in Volume 1 of Wiener's collected works (Wiener, 1976–1985, page 515), wrote as follows concerning Wiener's 1923 article:

> Having investigated the differential space from various directions, Wiener defines the Wiener measure as a $\sigma$-additive probability measure by means of Daniell's theory of integral.

It should not be thought, however, that Wiener defined a $\sigma$-additive probability measure and then found mean values as integrals with respect to that measure. Rather, as we just explained, he started with mean values and used Daniell's theory to obtain more. This Daniellian approach to probability, making mean value basic and probability secondary, has long taken a back seat to Kolmogorov's approach, but it still has its supporters (Haberman, 1996; Whittle, 2000).

### 3.3  Borel's Denumerable Probability

Impressive as it was and still is, Wiener's work played little role in the story leading to Kolmogorov's *Grundbegriffe*. The starring role was played instead by Borel.

In retrospect, Borel's use of measure theory in complex analysis in the 1890s already looks like probabilistic reasoning. Especially striking in this respect is the argument Borel gave for his claim that a Taylor series will usually diverge on the boundary of its circle of convergence (Borel, 1897). In general, he asserted, successive coefficients of the Taylor series, or at least successive groups of coefficients, are independent. He showed that each group of coefficients determines an arc on the circle, that the sum of lengths of the arcs diverges and that the Taylor series will diverge at a point on the circle if it belongs to infinitely many of the arcs. The arcs being independent and the sum of their lengths being infinite, a given point must be in infinitely many of them. To make sense of this argument, we must evidently take "in general" to mean that the coefficients are chosen at random and take "independent" to mean probabilistically independent; the conclusion then follows by what we now call the Borel–Cantelli lemma. Borel himself used probabilistic language when he reviewed this work in 1912 (Borel, 1912; Kahane, 1994). In the 1890s, however, Borel did not see complex analysis as a domain for probability, which is concerned with events in the real world.

In the new century, Borel did begin to explore the implications for probability of his and Lebesgue's work on measure and integration (Bru, 2001). His first comments came in an article in 1905 (Borel, 1905), where he pointed out that the new theory justified Poincaré's intuition that a point chosen at random from a line segment would be incommensurable with probability 1 and called attention to Anders Wiman's (1900, 1901) work on continued fractions, which had been inspired by the question of the stability of planetary motions, as an application of measure theory to probability.

Then, in 1909, Borel published a startling result—his strong law of large numbers (Borel, 1909a). This new result strengthened measure theory's connection both with geometric probability and with the heart of classical probability theory—the concept of independent trials. Considered as a statement in geometric probability, the law says that the fraction of 1's in the binary expansion of a real number chosen at random from $[0, 1]$ converges to $\frac{1}{2}$ with probability 1. Considered as a statement about independent trials (we may use the language of coin tossing, though Borel did not), it says that the fraction of heads in a denumerable sequence of independent tosses of a fair coin converges to $\frac{1}{2}$ with probability 1. Borel explained the geometric interpretation and he asserted that the result can be established using measure theory (Borel, 1909a, Section I.8). However, he set measure theory aside for philosophical reasons and provided an imperfect proof using denumerable versions of the rules of total and compound probability. It was left to others, most immediately Faber (1910, page 400) and Hausdorff (1914), to give rigorous measure-theoretic proofs (Doob, 1989, 1994; von Plato, 1994).

Borel's discomfort with a measure-theoretic treatment can be attributed to his unwillingness to assume countable additivity for probability (Barone and Novikoff, 1978; von Plato, 1994). He saw no logical absurdity in a countably infinite number of zero probabilities adding to a nonzero probability, and so instead of general appeals to countable additivity he preferred arguments that derive probabilities as limits as the number of trials increases (Borel, 1909a, Section I.4). Such arguments seemed to him stronger than formal appeals to countable additivity, because they exhibit the finitary pictures that are idealized by the infinitary pictures. He saw even more fundamental problems in the idea that Lebesgue measure can model a random choice (von Plato, 1994, pages 36–56; Knobloch, 2001). How can we choose a real number at random when most real numbers are not even definable in any constructive sense?

Although Hausdorff did not hesitate to equate Lebesgue measure with probability, his account of Borel's strong law, in his *Grundzüge der Mengenlehre* (Hausdorff, 1914, pages 419–421), treated it as a theorem about real numbers: the set of numbers in $[0, 1]$ with binary expansions for which the proportion of 1's converges to $\frac{1}{2}$ has Lebesgue measure 1. Later, Francesco Paolo Cantelli (1916a, b, 1917) rediscovered the strong law (he neglected, in any case, to cite Borel) and extended it to the more general result that the average of bounded random variables will converge to their mean with arbitrarily high probability. Cantelli's work inspired other authors to study the strong law and to sort out different concepts of probabilistic convergence.

By the early 1920s, it seemed to some that there were two different versions of Borel's strong law—one concerned with real numbers and one concerned with probability. Hugo Steinhaus (1923) proposed to clarify matters by axiomatizing Borel's theory of denumerable probability along the lines of Sierpiński's axiomatization of Lebesgue measure. Writing $A$ for the set of all infinite sequences of $\rho$'s and $\eta$'s ($\rho$ for "rouge" and $\eta$ for "noir"; now we are playing red or black rather than heads or tails), Steinhaus proposed the following axioms for a class $\mathfrak{K}$ of subsets of $A$ and a real-valued function $\mu$ that gives probabilities for the elements of $\mathfrak{K}$:

I. $\mu(E) \geq 0$ for all $E \in \mathfrak{K}$.
II. 1. For any finite sequence $e$ of $\rho$'s and $\eta$'s, the subset $E$ of $A$ consisting of all infinite sequences that begin with $e$ is in $\mathfrak{K}$.
    2. If two such sequences $e_1$ and $e_2$ differ in only one place, then $\mu(E_1) = \mu(E_2)$, where $E_1$ and $E_2$ are the corresponding sets.
    3. $\mu(A) = 1$.
III. $\mathfrak{K}$ is closed under finite and countable unions of disjoint elements, and $\mu$ is finitely and countably additive.
IV. If $E_1 \supset E_2$, and $E_1$ and $E_2$ are in $\mathfrak{K}$, then $E_1 \setminus E_2$ is in $\mathfrak{K}$.
V. If $E$ is in $\mathfrak{K}$ and $\mu(E) = 0$, then any subset of $E$ is in $\mathfrak{K}$.

Sierpiński's axioms for Lebesgue measure consisted of I, III, IV and V, together with an axiom that says that the measure $\mu(J)$ of an interval $J$ is its length. This last axiom being demonstrably equivalent to Steinhaus' axiom II, Steinhaus concluded that the theory of probability for an infinite sequence of binary trials is isomorphic with the theory of Lebesgue measure.

To show that his axiom II is equivalent to setting the measures of intervals equal to their length, Steinhaus used the Rademacher functions—the $n$th Rademacher function being the function that assigns a real number the value 1 or $-1$ depending on whether the $n$th digit in its dyadic expansion is 0 or 1. He also used these functions, which are independent random variables, in deriving Borel's strong law and related results. The work by Rademacher (1922) and Steinhaus marked the beginning of the Polish school of "independent functions," which made important contributions to probability theory during the period between the wars (Holgate, 1997).

### 3.4 Kolmogorov Enters the Stage

Although Steinhaus considered only binary trials in 1923, his reference to Borel's more general concept of denumerable probability pointed to generalizations. We find such a generalization in Kolmogorov's first article on probability, co-authored by Khinchin (Khinchin and Kolmogorov, 1925), which showed that a series of discrete random variables $y_1 + y_2 + \cdots$ will converge with probability 1 when the series of means and the series of variances both converge. The first section of the article, due to Khinchin, spells out how to represent the random variables as functions on $[0, 1]$: divide the interval into segments with lengths equal to the probabilities for $y_1$'s possible values, then divide each of these segments into smaller segments with lengths proportional to the probabilities for $y_2$'s possible values and so on. This, Khinchin noted with a nod to Rademacher and Steinhaus, reduces the problem to a problem about Lebesgue measure. This reduction was useful because the rules for working with Lebesgue measure were clear, while Borel's picture of denumerable probability remained murky.

Dissatisfaction with this detour into Lebesgue measure must have been one impetus for the *Grundbegriffe* (Doob, 1989, page 818). Kolmogorov made no such detour in his next article on the convergence of sums of independent random variables. In this sole-authored article (Kolmogorov, 1928), he took probabilities and expected values as his starting point, but even then he did not appeal to Fréchet's countably additive calculus. Instead, he worked with finite additivity and then stated an explicit ad hoc definition when he passed to a limit. For example, he defined the probability $P$ that the series $\sum_{n=1}^{\infty} y_n$ converges by the equation

$$P = \lim_{\eta \to 0} \lim_{n \to \infty} \lim_{N \to \infty} \mathfrak{W}\left[ \operatorname{Max} \left| \sum_{k=n}^{p} y_k \right|_{p=n}^{N} < \eta \right],$$

where $\mathfrak{W}(E)$ denotes the probability of the event $E$. [This formula does not appear in the Russian (Kolmogorov, 1986) and English (Kolmogorov, 1992) translations provided in Kolmogorov's collected works; there the argument has been modernized so as to eliminate it.] This recalls the way Borel proceeded in 1909: think through each passage to the limit.

It was in his seminal article on Markov processes (Kolmogorov, 1931) that Kolmogorov first explicitly and freely used Fréchet's calculus as his framework for probability. In this article, Kolmogorov considered a system with a set of states $\mathfrak{A}$. For any two time points $t_1$ and $t_2$ ($t_1 < t_2$), any state $x \in \mathfrak{A}$ and any element $\mathfrak{E}$ in a collection $\mathfrak{F}$ of subsets of $\mathfrak{A}$, he wrote $P(t_1, x, t_2, \mathfrak{E})$ for the probability, when the system is in state $x$ at time $t_1$, that it will be in a state in $\mathfrak{E}$ at time $t_2$. Citing Fréchet, Kolmogorov assumed that $P$ is countably additive as a function of $\mathfrak{E}$ and that $\mathfrak{F}$ is closed under differences and countable unions, and contains the empty set, all singletons and $\mathfrak{A}$. However, the focus was not on Fréchet; it was on the equation that ties together the transition probabilities, now called the Chapman–Kolmogorov equation. The article launched the study of this equation by purely analytical methods, a study that kept probabilists occupied for 50 years.

As many commentators have noted, the 1931 article makes no reference to probabilities for trajectories. There is no suggestion that such probabilities are needed for a stochastic process to be well defined. Consistent transition probabilities, it seems, are enough. Bachelier (1900, 1910, 1912) is cited as the first to study continuous-time stochastic processes, but Wiener is not cited.

## 4. HILBERT'S SIXTH PROBLEM

At the beginning of the twentieth century, many mathematicians were dissatisfied with what they saw as a lack of clarity and rigor in the probability calculus. The whole calculus seemed to be concerned with concepts that lie outside mathematics: event, trial, randomness, probability. As Henri Poincaré wrote, "one can hardly give a satisfactory definition of probability" (Poincaré, 1912, page 24).

The most celebrated call for clarification came from David Hilbert. The sixth of the twenty-three open problems that Hilbert presented to the International Congress of Mathematicians in Paris in 1900 was to treat axiomatically, after the model of geometry, those parts of physics in which mathematics already played an outstanding role, especially probability and mechanics (Hilbert, 1902; Hochkirchen, 1999). To explain what he meant by axioms for probability, Hilbert cited Georg Bohlmann, who had labeled the rules of total and compound probability axioms rather than theorems in his lectures on the mathematics of life insurance (Bohlmann, 1901). In addition to a logical investigation of these axioms, Hilbert called for a "rigorous and satisfactory development of the method of average values in mathematical physics, especially in the kinetic theory of gases."

Hilbert's call for a mathematical treatment of average values was answered in part by the work on integration that we discussed in the preceding section, but his suggestion that the classical rules for probability should be treated as axioms on the model of geometry was an additional challenge. Among the early responses, we may mention the following:

- In his Zürich dissertation, Rudolf Laemmel (1904) discussed the rules of total and compound probability as axioms, but he stated the rule of compound probability only in the case of independence, a concept he did not explicate. (For excerpts, see Schneider, 1988, pages 359–366.)

- In his Göttingen dissertation, directed by Hilbert himself, Ugo Broggi (1907) gave only two axioms: an axiom stating that the sure event has probability 1, and an axiom stating the rule of total probability. Following tradition, he then defined probability as a ratio (a ratio of numbers of cases in the discrete setting; a ratio of the Lebesgue measures of two sets in the geometric setting) and verified his axioms. He did not state an axiom that corresponds to the classical rule of compound probability. Instead, he gave this name to a rule for calculating the probability of a Cartesian product, which he derived from the definition of geometric probability in terms of Lebesgue measure. (For excerpts, see Schneider, 1988, pages 367–377.) Broggi mistakenly claimed that his axiom of total probability (finite additivity) implied countable additivity (Steinhaus, 1923).

- In an article written in 1920, published in 1923 and listed in the bibliography of the *Grundbegriffe*, Antoni Łomnicki (1923) proposed that probability should always be understood relative to a density $\phi$ on a set $\mathcal{M}$ in $R^r$. Łomnicki defined this probability by combining two of Carathéodory's ideas: the idea of $p$-dimensional measure and the idea of defining the integral of a function on a set as the measure of the region between the set and the function's graph (see Section 3.1 above). The probability of a subset $m$ of $\mathcal{M}$, according to Łomnicki, is the ratio of the measure of the region between $m$ and $\phi$'s graph to the measure of the region between $\mathcal{M}$ and this graph. If $\mathcal{M}$ is an $r$-dimensional subset of $R^r$, then the measure being used is Lebesgue measure on $R^{r+1}$; if $\mathcal{M}$ is a lower dimensional subset of $R^r$, say $p$-dimensional, then the measure is the $(p+1)$-dimensional Carathéodory measure. This definition covers discrete as well as continuous probability: in the discrete case, $\mathcal{M}$ is a set of discrete points, the function $\phi$ assigns each point its probability, and the region between a subset $m$ and the graph of $\phi$ consists of a line segment for each point in $m$, whose Carathéodory measure is its length (i.e., the point's probability). The rule of total probability follows. Like Broggi, Łomnicki treated the rule of compound probability as a rule for relating probabilities on a Cartesian product to probabilities on its components. He did not consider it an axiom, because it holds only if the density itself is a product density.

• In an article published in Russian, Sergei Bernstein (1917) showed that probability theory can be founded on qualitative axioms for numerical coefficients that measure the probabilities of propositions. He also developed this idea in his probability textbook (Bernstein, 1927), and Kolmogorov listed both the article and the book in the bibliography of the *Grundbegriffe*. John Maynard Keynes included Bernstein's article in the bibliography of his probability book (Keynes, 1921), but Bernstein's work was subsequently ignored by English-language authors on qualitative probability. It was first summarized in English in Samuel Kotz's translation of Leonid E. Maistrov's (1974) history of probability.

We now discuss at greater length responses by von Mises, Slutsky, Kolmogorov and Cantelli.

### 4.1 von Mises' Collectives

The concept of a collective was introduced into the German scientific literature by Gustav Fechner's (1897) *Kollektivmasslehre*, which appeared ten years after the author's death. The concept was quickly taken up by Georg Helm (1902) and Heinrich Bruns (1906).

Fechner wrote about the concept of a *Kollektivgegenstand* (collective object) or a *Kollektivreihe* (collective series). It was only later, in Meinong (1915) for example, that we see these names abbreviated to *Kollektiv*. As the name *Kollektivreihe* indicates, a *Kollektiv* is a population of individuals given in a certain order; Fechner called the ordering the *Urliste*. It was supposed to be irregular—random, we would say. Fechner was a practical scientist, not concerned with the theoretical notion of probability, but as Helm and Bruns realized, probability theory provides a framework for studying collectives.

The concept of a collective was developed by Richard von Mises (1919, 1928, 1931). His contribution was to realize that the concept can be made into a mathematical foundation for probability theory. As von Mises defined it, a collective is an infinite sequence of outcomes that satisfies two axioms:

1. The relative frequency of each outcome converges to a real number (the probability of the outcome) as we look at longer and longer initial segments of the sequence.

2. The relative frequency converges to the same probability in any subsequence selected without knowledge of the future (we may use knowledge of the outcomes so far in deciding whether to include the next outcome in the subsequence).

The second property says we cannot change the odds by selecting a subsequence of trials on which to bet; this is von Mises' version of the "hypothesis of the impossibility of a gambling system," and it assures the irregularity of the *Urliste*.

According to von Mises, the purpose of the probability calculus is to identify situations where collectives exist and the probabilities in them are known, and to derive probabilities for other collectives from these given probabilities. He pointed to three domains where probabilities for collectives are known: (1) games of chance where devices are carefully constructed so the axioms will be satisfied, (2) statistical phenomena where the two axioms can be confirmed, to a reasonable degree and (3) branches of theoretical physics where the two axioms play the same hypothetical role as other theoretical assumptions (von Mises, 1931, pages 25–27).

von Mises derived the classical rules of probability, such as the rules for adding and multiplying probabilities, from rules for constructing new collectives from an initial one. He had several laws of large numbers. The simplest was his definition of probability: the probability of an event is the event's limiting frequency in a collective. Others arose as one constructed further collectives.

The ideas of von Mises were taken up by a number of mathematicians in the 1920s and 1930s. Kolmogorov's bibliography includes an article by Arthur Copeland (1932) that proposed founding probability theory on particular rules for selecting subsequences in von Mises' scheme, as well as articles by Karl Dörge (1930), Hans Reichenbach (1932) and Erhard Tornier (1933) that argued for related schemes. But the most prominent mathematicians of the time, including the Göttingen mathematicians (Mac Lane, 1995), the French probabilists and the British statisticians, were hostile or indifferent.

Collectives were given a rigorous mathematical basis by Abraham Wald (1938) and Alonzo Church (1940), but the claim that they provide a foundation for probability was refuted by Jean Ville (1939). Ville pointed out that whereas a collective in von Mises' sense will not be vulnerable to a gambling system that chooses a subsequence of trials on which to bet, it may still be vulnerable to a more clever gambling system, which also varies the amount of the bet and the outcome on which to bet.

## 4.2 Slutsky's Calculus of Valences and Kolmogorov's General Theory of Measure

In an article published in Russian Evgeny Slutsky (1922) presented a viewpoint that greatly influenced Kolmogorov. As Kolmogorov (1948) said in an obituary for Slutsky, Slutsky was "the first to give the right picture of the purely mathematical content of probability theory."

How do we make probability purely mathematical? Markov had claimed to do this in his textbook, but Slutsky did not think Markov had succeeded, because Markov had retained the subjective notion of equipossibility. The solution, Slutsky felt, was to remove both the word "probability" and the notion of equally likely cases from the theory. Instead of beginning with equally likely cases, one should begin by assuming merely that numbers are assigned to cases and that when a case assigned the number $\alpha$ is further subdivided, the numbers assigned to the subcases should add to $\alpha$. The numbers assigned to cases might be equal or they might not. The addition and multiplication theorems would be theorems in this abstract calculus, but it should not be called the probability calculus. In place of "probability," he suggested the unfamiliar word валентность, or "valence." (Laemmel had earlier used the German *valenz*.) Probability would be only one interpretation of the calculus of valences, a calculus fully as abstract as group theory.

Slutsky listed three distinct interpretations of the calculus of valences:

1. Classical probability (equally likely cases).
2. Finite empirical sequences (frequencies).
3. Limits of relative frequencies. (Slutsky remarked that this interpretation is particularly popular with the English school.)

Slutsky did not think probability could be reduced to limiting frequency, because sequences of independent trials have properties that go beyond their possessing limiting frequencies. Initial segments of the sequences have properties that are not imposed by the eventual convergence of the frequency, and the sequences must be irregular in a way that resists the kind of selection discussed by von Mises.

Slutsky's idea that probability could be an instance of a broader abstract theory was taken up by Kolmogorov in a thought piece in Russian (Kolmogorov, 1929), before his forthright use of Fréchet's theory in his article on Markov processes in 1930 (Kolmogorov, 1931). Whereas Slutsky had mentioned frequencies as an alternative interpretation of a general calculus, Kolmogorov pointed to more mathematical examples: the distribution of digits in the decimal expansions of irrationals, Lebesgue measure in an $n$-dimensional cube and the density of a set $A$ of positive integers (the limit as $n \to \infty$ of the fraction of the integers between 1 and $n$ that are in $A$).

The abstract theory Kolmogorov sketches is concerned with a function $M$ that assigns a nonnegative number $M(E)$ to each element $E$ of a class of subsets of a set $A$. He called $M(E)$ the measure (мера) of $E$ and he called $M$ a measure specification (мероопределение). So as to accommodate all the mathematical examples he had in mind, he assumed, in general, neither that $M$ is countably additive nor that the class of subsets to which it assigns numbers is a field. Instead, he assumed only that when $E_1$ and $E_2$ are disjoint and $M$ assigns a number to two of the three sets $E_1$, $E_2$ and $E_1 \cup E_2$, it also assigns a number to the third, and that

$$M(E_1 \cup E_2) = M(E_1) + M(E_2)$$

then holds (cf. Steinhaus' axioms III and IV). In the case of probability, however, he did suggest (using different words) that $M$ should be countably additive and that the class of subsets to which it assigns numbers should be a field, for only then can we uniquely define probabilities for countable unions and intersections, and this seems necessary to justify arguments involving events such as the convergence of random variables.

He defined the abstract Lebesgue integral of a function $f$ on $A$, and he commented that countable additivity is to be assumed whenever such an integral is discussed. He wrote $M_{E_1}(E_2) = M(E_1 E_2)/M(E_1)$ "by analogy with the usual concept of relative probability." He defined independence for partitions, and he commented, no doubt in reference to Borel's strong law and other results in number theory, that the notion of independence is responsible for the power of probabilistic methods within pure mathematics.

The mathematical core of the *Grundbegriffe* is already here. Many years later, in his commentary in Volume II of his collected works (Kolmogorov, 1992, page 520), Kolmogorov said that only the set-theoretic treatment of conditional probability and the theory of distributions in infinite products were missing. Also missing, though, is the bold rhetorical move that Kolmogorov made in the *Grundbegriffe*—giving the abstract theory the name probability.

### 4.3 The Axioms of Steinhaus and Ulam

In the 1920s and 1930s, the city of Lwów in Poland was a vigorous center of mathematical research, led by Hugo Steinhaus. (Though it was in Poland between the two World Wars, Lwów is now in Ukraine. Its name is spelled differently in different languages: Lwów in Polish, Lviv in Ukrainian and Lvov in Russian. When part of Austria–Hungary and, briefly, Germany, it was Lemberg. Some articles in our bibliography refer to it as Léopol.) In 1929, Steinhaus' work on limit theorems intersected with Kolmogorov's, and his approach promoted the idea that probability should be axiomatized in the style of measure theory.

As we saw in Section 3.3, Steinhaus had already, in 1923, formulated axioms for heads and tails isomorphic to Sierpiński's axioms for Lebesgue measure. This isomorphism had more than a philosophical purpose; Steinhaus used it to prove Borel's strong law. In a pair of articles written in 1929 and published in 1930 (Steinhaus, 1930a, b), Steinhaus extended his approach to limit theorems that involved an infinite sequence of independent draws $\theta_1, \theta_2, \ldots$ from the interval $[0, 1]$. His axioms for this case were the same as for the binary case (Steinhaus, 1930b, pages 22–23), except that the second axiom, which determines probabilities for initial finite sequences of heads and tails, was replaced by an axiom that determines probabilities for initial finite sequences $\theta_1, \theta_2, \ldots, \theta_n$:

> The probability that $\theta_i \in \Theta_i$ for $i = 1, \ldots, n$, where the $\Theta_i$ are measurable subsets of $[0, 1]$, is
>
> $$|\Theta_1| \cdot |\Theta_2| \cdots |\Theta_n|,$$
>
> where $|\Theta_i|$ is the Lebesgue measure of $\Theta_i$.

Steinhaus presented his axioms as a "logical extrapolation" of the classical axioms to the case of an infinite number of trials (Steinhaus, 1930b, page 23). They were more or less tacitly used, he asserted, in all classical problems, such as the problem of the gambler's ruin, where the game as a whole—not merely finitely many rounds—must be considered (Steinhaus, 1930a, page 409). As in the case of heads and tails, Steinhaus showed that there are probabilities that uniquely satisfy his axioms by setting up an isomorphism with Lebesgue measure on $[0, 1]$, this time using a sort of Peano curve to map $[0, 1]^\infty$ onto $[0, 1]$. He used the isomorphism to prove several limit theorems, including one that formalized Borel's 1897 claim concerning the circle of convergence of a Taylor series with randomly chosen coefficients.

Steinhaus' axioms were measure-theoretic, but they were not yet abstract. His words suggested that his ideas should apply to all sequences of random variables, not merely ones uniformly distributed, and he even considered the case where the variables were complex-valued rather than real-valued, but he did not step outside the geometric context to consider probability on abstract spaces. This step was taken by Stanisław Ulam, one of Steinhaus' junior colleagues at Lwów. At the International Congress of Mathematicians in Zürich in 1932, Ulam announced that he and another Lwów mathematician, Zbigniew Łomnicki (a nephew of Antoni Łomnicki), had shown that product measures can be constructed in abstract spaces (Ulam, 1932).

Ulam and Łomnicki's axioms for a measure $m$ were simple. We can put them in today's language by saying that $m$ is a probability measure on a $\sigma$-algebra that is complete (includes all null sets) and contains all singletons. Ulam announced that from a countable sequence of spaces with such probability measures, one can construct a probability measure that satisfies the same conditions on the product space.

We do not know whether Kolmogorov knew about Ulam's announcement when he wrote the *Grundbegriffe*. Ulam's axioms would have held no novelty for him, but he would presumably have found the result on product measures interesting. When it finally appeared, Łomnicki and Ulam (1934) listed the same axioms as Ulam's announcement had, but it now cited the *Grundbegriffe* as authority for them. Kolmogorov (1935) cited their article in turn in a short list of introductory literature in mathematical probability.

### 4.4 Cantelli's Abstract Theory

Like Borel, Castelnuovo and Fréchet, Francesco Paolo Cantelli turned to probability after distinguishing himself in other areas of mathematics. It was only in the 1930s, about the same time as the *Grundbegriffe* appeared, that he introduced his own abstract theory of probability. This theory, which has important affinities with Kolmogorov's, is developed most clearly in an article included in the *Grundbegriffe*'s bibliography (Cantelli, 1932) and a lecture he gave in 1933 at the Institut Henri Poincaré in Paris (Cantelli, 1935).

Cantelli (1932) argued for a theory that makes no appeal to empirical notions such as possibility, event, probability or independence. This abstract theory, he said, should begin with a set of points that have finite nonzero measure. This could be any set for which measure is defined, perhaps a set of points on a surface. He wrote $m(E)$ for the area of a subset $E$. He noted that $m(E_1 \cup E_2) = m(E_1) + m(E_2)$, provided $E_1$ and $E_2$ are disjoint, and $0 \le m(E_1E_2)/m(E_i) \le 1$ for $i = 1, 2$. He called $E_1$ and $E_2$ *multipliable* when $m(E_1E_2) = m(E_1)m(E_2)$. Much of probability theory, he noted, including Bernoulli's law of large numbers and Khinchin's law of the iterated logarithm, can be carried out at this abstract level.

Cantelli (1935) explained how his abstract theory relates to frequencies in the world. The classical calculus of probability, he said, should be developed for a particular class of events in the world in three steps:

1. Study experimentally the equally likely cases (check that they happen equally frequently), thus justifying experimentally the rules of total and compound probability.
2. Develop an abstract theory based only on the rules of total and compound probability, without reference to their empirical justification.
3. Deduce probabilities from the abstract theory and use them to predict frequencies.

His own theory, Cantelli explains, is the one obtained in the second step.

Cantelli's 1932 article and 1933 lecture were not really sources for the *Grundbegriffe*. Kolmogorov's earlier work (Kolmogorov, 1929, 1931) had already went well beyond anything Cantelli did in 1932, in both degree of abstraction and mathematical clarity. The 1933 lecture was more abstract, but obviously came too late to influence the *Grundbegriffe*. However, Cantelli did develop independently of Kolmogorov the project of combining a frequentist interpretation of probability with an abstract axiomatization that retained in some form the classical rules of total and compound probability. This project had been in the air for 30 years.

## 5. THE GRUNDBEGRIFFE

The *Grundbegriffe* was an exposition, not another research contribution. In his preface, after acknowledging Fréchet's work, Kolmogorov said this:

> In the pertinent mathematical circles it has been common for some time to construct probability theory in accordance with this general point of view. But a complete presentation of the whole system, free from superfluous complications, has been missing (though a book by Fréchet, [2] in the bibliography, is in preparation).

Kolmogorov aimed to fill this gap, and he did so brilliantly and concisely, in just 62 pages. Fréchet's much longer book, which finally appeared in two volumes (Fréchet, 1937–1938), is regarded by some as a mere footnote to Kolmogorov's achievement.

Fréchet's own evaluation of the *Grundbegriffe*'s contribution, quoted at the beginning of this article, is correct so far as it goes. Borel had introduced countable additivity into probability in 1909, and in the following 20 years, many authors, including Kolmogorov, had explored its consequences. The *Grundbegriffe* merely rounded out the picture by explaining that nothing more was needed. However, Kolmogorov's mathematical achievement, especially his definitive work on the classical limit theorems, had given him the grounds and the authority to say that nothing more was needed.

Moreover, Kolmogorov's appropriation of the name *probability* was an important rhetorical achievement, with enduring implications. Slutsky in 1922 and Kolmogorov himself in 1927 had proposed a general theory of additive set functions but had relied on the classical theory to say that probability should be a special case of this general theory. Now Kolmogorov proposed axioms for probability. The numbers in his abstract theory were probabilities, not merely valences or меры. His philosophical justification for proceeding in this way so resembled the justification that Borel and Lévy had given for the classical theory that they could hardly take exception.

It was not really true that nothing more was needed. Those who studied Kolmogorov's formulation in detail soon realized that his axioms and definitions were inadequate in a number of ways. Most saliently, his treatment of conditional probability was not adequate for the burgeoning theory of Markov processes. In addition, there were other points in the monograph where he could not obtain natural results at the abstract level and had to fall back to the classical examples—discrete probabilities and probabilities in Euclidean spaces. These shortcomings only gave impetus to the new theory, because the project of filling in the gaps provided exciting work for a new generation of probabilists.

In this section we take a fresh look at the *Grundbegriffe*. We review its six axioms and two ideas that were, as Kolmogorov himself pointed out in his preface, novel at the time: the construction of probabilities on infinite-dimensional spaces (his famous consistency theorem) and the definition of conditional probability using the Radon–Nikodym theorem. Then we look at the explicitly philosophical part of the monograph: the two pages in Chapter I where Kolmogorov explains the empirical origin and meaning of his axioms.

## 5.1 The Mathematical Framework

Kolmogorov's six axioms for probability are so familiar that it seems superfluous to repeat them, but so concise that it is easy to do so. We do repeat them and then we discuss the two points just mentioned: the consistency theorem and the treatment of conditional probability and expectation. As we will see, the mathematics was due to earlier authors—Daniell in the case of the consistency theorem and Nikodym in the case of conditional probabilities and expectations. Kolmogorov's contribution, more rhetorical and philosophical than mathematical, was to bring this mathematics into a framework for probability.

5.1.1 *The six axioms.* Kolmogorov began with five axioms concerning a set $E$ and a set $\mathfrak{F}$ of subsets of $E$, which he called *random events*:

   I. $\mathfrak{F}$ is a field of sets.
  II. $\mathfrak{F}$ contains the set $E$.

III. To each set $A$ from $\mathfrak{F}$ is assigned a nonnegative real number $\mathsf{P}(A)$. This number $\mathsf{P}(A)$ is called the probability of the event $A$.

IV. The $\mathsf{P}(E) = 1$.

V. If $A$ and $B$ are disjoint, then

$$\mathsf{P}(A \cup B) = \mathsf{P}(A) + \mathsf{P}(B).$$

He then added a sixth axiom, redundant for finite $\mathfrak{F}$ but independent of the first five axioms for infinite $\mathfrak{F}$:

VI. If $A_1 \supseteq A_2 \supseteq \cdots$ is a decreasing sequence of events from $\mathfrak{F}$ with $\bigcap_{n=1}^{\infty} A_n = \varnothing$, then $\lim_{n \to \infty} \mathsf{P}(A_n) = 0$.

This is the *axiom of continuity*. Given the first five axioms, it is equivalent to countable additivity.

The six axioms can be summarized by saying that $\mathsf{P}$ is a nonnegative additive set function in the sense of Fréchet with $\mathsf{P}(E) = 1$.

Unlike Fréchet, who had debated countable additivity with de Finetti (Fréchet, 1930; de Finetti, 1930; Cifarelli and Regazzini, 1996), Kolmogorov did not make a substantive argument for it. Instead, he said this (page 14):

> . . . Since the new axiom is essential only for infinite fields of probability, it is hardly possible to explain its empirical meaning. . . . In describing any actual observable random process, we can obtain only finite fields of probability. Infinite fields of probability occur only as idealized models of real random processes. *This understood, we limit ourselves arbitrarily to models that satisfy Axiom VI.* So far this limitation has been found expedient in the most diverse investigations.

This echoes Borel who adopted countable additivity not as a matter of principle but because he had not encountered circumstances where its rejection seemed expedient (Borel, 1909a, Section I.5). However, Kolmogorov articulated even more clearly than Borel the purely instrumental significance of infinity.

5.1.2 *Probability distributions in infinite-dimensional spaces.* Suppose, using modern terminology, that $(E_1, \mathfrak{F}_1), (E_2, \mathfrak{F}_2), \ldots$ is a sequence of measurable spaces. For each finite set of indices, say $i_1, \ldots, i_n$, write $\mathfrak{F}^{i_1, \ldots, i_n}$ for the induced $\sigma$-algebra in the product space $\prod_{j=1}^{n} E_{i_j}$. Write $E$ for the product of all the $E_i$ and write $\mathfrak{F}$ for the algebra (not a $\sigma$-algebra) that consists of all the cylinder subsets of $E$ corresponding to elements of the various $\mathfrak{F}^{i_1, \ldots, i_n}$. Suppose we define consistent probability measures for all the marginal spaces $(\prod_{j=1}^{n} E_{i_j}, \mathfrak{F}^{i_1, \ldots, i_n})$. This defines a set function on $(E, \mathfrak{F})$. Is it countably additive?

In general, the answer is negative; a counterexample was given by Erik Sparre Andersen and Børge Jessen in 1948, but as we noted in Section 4.3, Ulam had given a positive answer for the case where the marginal measures are product measures. Kolmogorov's consistency theorem, in Section 4 of Chapter III of the *Grundbegriffe*, gave a positive answer for another case, where each $E_i$ is a copy of the real numbers and each $\mathfrak{F}_i$ consists of the Borel sets. (Formally, we should acknowledge, Kolmogorov had a slightly different starting point: finite-dimensional distribution functions, not finite-dimensional measures.)

In his September 1919 article (Daniell, 1919b), Daniell had proven a closely related theorem. Although Kolmogorov did not cite Daniell in the *Grundbegriffe*, the essential mathematical content of Kolmogorov's result is already in Daniell's. This point was recognized quickly; Jessen (1935) called attention to Daniell's priority in an article that appeared in MIT's *Journal of Mathematics and Physics*, together with an article by Wiener that also called attention to Daniell's result. In a commemoration of Kolmogorov's early work, Doob (1989) hazards the guess that Kolmogorov was unaware of Daniell's result when he wrote the *Grundbegriffe*. This may be true. He would not have been the first author to repeat Daniell's work; Jessen had presented the result as his own to the Seventh Scandinavian Mathematical Conference in 1929 and had become aware of Daniell's priority only in time to acknowledge it in a footnote to his contribution to the proceedings (Jessen, 1930).

It is implausible that Kolmogorov was still unaware of Daniell's construction after the comments by Wiener and Jessen, but in 1948 he again ignored Daniell while claiming the construction of probability measures on infinite products as a Soviet achievement (Gnedenko and Kolmogorov, 1948, Section 3.1). Perhaps this can be dismissed as mere propaganda, but we should also remember that the *Grundbegriffe* was not meant as a contribution to pure mathematics. Daniell's and Kolmogorov's theorems seem almost identical when they are assessed as mathematical discoveries, but they differed in context and purpose. Daniell was not thinking about probability, whereas the slightly different theorem formulated by Kolmogorov was about probability. Neither Daniell nor Wiener undertook to make probability into a conceptually independent branch of mathematics by establishing a general method for representing it measure-theoretically.

Kolmogorov's theorem was more general than Daniell's in one respect—Kolmogorov considered an index set of arbitrary cardinality, whereas Daniell considered only denumerable cardinality. This greater generality is merely formal, in two senses: it involves no additional mathematical complications and it has no practical use. The obvious use of a nondenumerable index would be to represent continuous time, and so we might conjecture that Kolmogorov was thinking of making probability statements about trajectories, as Wiener had done in the 1920s. However, Kolmogorov's construction does not accomplish anything in this direction. The $\sigma$-algebra on the product obtained by the construction contains too few sets; in the case of Brownian motion, it does not include the set of continuous trajectories. It took some decades of further research to develop general methods of extension to $\sigma$-algebras rich enough to include the infinitary events one typically wants to discuss (Doob, 1953; Bourbaki, 1994, pages 243–245). The topological character of these extensions and the failure of the consistency theorem for arbitrary Cartesian products remain two important caveats to the *Grundbegriffe*'s thesis that probability is adequately represented by the abstract notion of a probability measure.

### 5.1.3 *Experiments and conditional probability.*

In the case where $A$ has nonzero probability, Kolmogorov defined $\mathsf{P}_A(B)$ in the usual way. He called it *bedingte Wahrscheinlichkeit*, which translates into English as "conditional probability." Before the *Grundbegriffe*, this term was less common than "relative probability."

Kolmogorov's treatment of conditional probability and expectation was novel. It began with a set-theoretic formalization of the concept of an *experiment* (*Versuch* in German). Here Kolmogorov had in mind a subexperiment of the grand experiment defined by the conditions $\mathfrak{S}$. The subexperiment may give only limited information about the outcome $\xi$ of the grand experiment. It defines a partition $\mathfrak{A}$ of the sample space $E$ for the grand experiment: its outcome amounts to specifying which element of $\mathfrak{A}$ contains $\xi$. Kolmogorov formally identified the subexperiment with $\mathfrak{A}$. Then he introduced the idea of conditional probability relative to $\mathfrak{A}$:

- In the finite case, he wrote $\mathsf{P}_{\mathfrak{A}}(B)$ for the random variable whose value at each point $\xi$ of $E$ is $\mathsf{P}_A(B)$, where $A$ is the element of $\mathfrak{A}$ that contains $\xi$, and he called this random variable the "conditional probability of $B$ after the experiment $\mathfrak{A}$" (page 12). This random variable is well defined for all the $\xi$ in elements of $\mathfrak{A}$ that have positive probability, and these $\xi$ form an event that has probability 1.

- In the general case, he represented the partition $\mathfrak{A}$ by a function $u$ on $E$ that induces it and he wrote $\mathsf{P}_u(B)$ for any random variable that satisfies

$$\mathsf{P}_{\{u \subset A\}}(B) = \mathsf{E}_{\{u \subset A\}}\mathsf{P}_u(B)$$

for every set $A$ of possible values of $u$ such that the subset $\{\xi | u(\xi) \in A\}$ of $E$ (this is what he meant by $\{u \subset A\}$) is measurable and has positive probability (page 42). By the Radon–Nikodym theorem (only recently proven by Nikodym), this random variable is unique up to a set of probability 0. Kolmogorov called it the "conditional probability of $B$ with respect to (or knowing) $u$." He defined $\mathsf{E}_u(y)$, which he called "the conditional expectation of the variable $y$ for a known value of $u$," analogously (page 46).

Kolmogorov was doing no new mathematics here; the mathematics is Nikodym's. However, Kolmogorov was the first to point out that Nikodym's result can be used to derive conditional probabilities from absolute probabilities.

We should not, incidentally, jump to the conclusion that Kolmogorov had abandoned the emphasis on transition probabilities he had displayed in his 1931 article and now wanted to start the study of stochastic processes with unconditional probabilities. Even in 1935, he recommended the opposite (Kolmogorov, 1935, pages 168–169 of the English translation).

5.1.4 *When is conditional probability meaningful?* To illustrate his understanding of conditional probability, Kolmogorov discussed Bertrand's paradox of the great circle, which he called, with no specific reference, a Borelian paradox. His explanation of the paradox was simple but formal. After noting that the probability distribution for the second point conditional on a particular great circle is not uniform, he said:

> This demonstrates the inadmissibility of the idea of conditional probability with respect to a given isolated hypothesis with probability zero. One obtains a probability distribution for the latitude on a given great circle only when that great circle is considered as an element of a partition of the entire surface of the sphere into great circles with the given poles (page 45).

This explanation has become part of the culture of probability theory, but it cannot completely replace the more substantive explanations given by Borel.

Borel insisted that we explain how the measurement on which we will condition is to be carried out. This accords with Kolmogorov's insistence that a partition be specified, because a procedure for measurement will determine such a partition. Kolmogorov's explicitness on this point was a philosophical advance. On the other hand, Borel demanded more than the specification of a partition. He demanded that the measurement be specified realistically enough that we can see partitions into events of positive probability, not just a theoretical limiting partition into events of probability 0.

Borel's demand that we be told how the theoretical partition into events of probability 0 arises as a limit of partitions into events of positive probability again compromises the abstract picture by introducing topological ideas, but this seems to be needed so as to rule out nonsense. This point was widely discussed in the 1940s and 1950s. Dieudonné (1948) and Lévy (1959) gave examples in which the conditional probabilities defined by Kolmogorov do not have versions (functions of $\xi$ for fixed $B$) that form sensible probability measures (when considered as functions of $B$ for fixed $\xi$). Gnedenko and Kolmogorov (1949) and Blackwell (1956) formulated conditions on measurable spaces or probability measures that rule out such examples. For modern formulations of these conditions, see Rogers and Williams (2000).

## 5.2  The empirical origin of the axioms

Kolmogorov devoted about two pages of the *Grundbegriffe* to the relation between his axioms and the real world. These two pages, a concise statement of Kolmogorov's frequentist philosophy, are so important to our story that we quote them in full. We then discuss how this philosophy was related to the thinking of his predecessors and how it fared in the decades following 1933.

5.2.1  *In Kolmogorov's own words.* Section 2 of Chapter I of the *Grundbegriffe* is titled "*Das Verhältnis zur Erfahrungswelt.*" It is only two pages in length. This subsection consists of a translation of the section in its entirety.

### The relation to the world of experience

The theory of probability is applied to the real world of experience as follows:

1.  Suppose we have a certain system of conditions $\mathfrak{S}$, capable of unlimited repetition.

2.  We study a fixed circle of phenomena that can arise when the conditions $\mathfrak{S}$ are realized. In general, these phenomena can come out in different ways in different cases where the conditions are realized. Let $E$ be the set of the different possible variants $\xi_1, \xi_2, \ldots$ of the outcomes of the phenomena. Some of these variants might actually not occur. We include in the set $E$ all the variants we regard *a priori* as possible.

3.  If the variant that actually appears when conditions $\mathfrak{S}$ are realized belongs to a set $A$ that we define in some way, then we say that the event $A$ has taken place.

EXAMPLE.  The system of conditions $\mathfrak{S}$ consists of flipping a coin twice. The circle of phenomena mentioned in point 2 consists of the appearance, on each flip, of heads or tails. It follows that there are four possible variants (*elementary events*), namely

<div align="center">heads—heads, heads—tails,</div>

<div align="center">tails—heads, tails—tails.</div>

Consider the event $A$ that there is a repetition. This event consists of the first and fourth elementary events. Every event can similarly be regarded as a set of elementary events.

4.  Under certain conditions, that we will not go into further here, we may assume that an event $A$ that does or does not occur under conditions $\mathfrak{S}$ is assigned a real number $\mathsf{P}(A)$ with the following properties:

A. One can be practically certain that if the system of conditions $\mathfrak{S}$ is repeated a large number of times, $n$, and the event $A$ occurs $m$ times, then the ratio $m/n$ will differ only slightly from $\mathsf{P}(A)$.

B. If $\mathsf{P}(A)$ is very small, then one can be practically certain that the event $A$ will not occur on a single realization of the conditions $\mathfrak{S}$.

*Empirical deduction of the axioms.* Usually one can assume that the system $\mathfrak{F}$ of events $A, B, C \ldots$ that come into consideration and are assigned definite probabilities forms a field that contains $E$ (Axioms I and II and the first half of Axiom III—the existence of the probabilities). It is further evident that $0 \le m/n \le 1$ always holds, so that the second half of Axiom III appears completely natural. We always have $m = n$ for the event $E$, so we naturally set $\mathsf{P}(E) = 1$ (Axiom IV). Finally, if $A$ and $B$ are mutually incompatible (in other words, the sets $A$ and $B$ are disjoint), then $m = m_1 + m_2$, where $m$, $m_1$ and $m_2$ are the numbers of experiments in which the events $A \cup B$, $A$ and $B$ happen, respectively. It follows that

$$\frac{m}{n} = \frac{m_1}{n} + \frac{m_2}{n}.$$

So it appears appropriate to set $\mathsf{P}(A \cup B) = \mathsf{P}(A) + \mathsf{P}(B)$.

Remark I.   If two assertions are both practically certain, then the assertion that they are simultaneously correct is practically certain, though with a little lower degree of certainty. But if the number of assertions is very large, we cannot draw any conclusion whatsoever about making the assertions simultaneously from the practical certainty of each of them individually. So it in no way follows from Principle A that $m/n$ will differ only a little from $\mathsf{P}(A)$ in every one of a very large number of series of experiments, where each series consists of $n$ experiments.

Remark II.   By our axioms, the impossible event (the empty set) has the probability $\mathsf{P}(\varnothing) = 0$. But the converse inference, from $\mathsf{P}(A) = 0$ to the impossibility of $A$, does not by any means follow. By Principle B, the event $A$'s having probability zero implies only that it is practically impossible that it will happen on a particular unrepeated realization of the conditions $\mathfrak{S}$. This by no means implies that the event $A$ will not appear in the course of a sufficiently long series of experiments. When $\mathsf{P}(A) = 0$ and $n$ is very large, we can only say, by Principle A, that the quotient $m/n$ will be very small—it might, for example, be equal to $1/n$.

5.2.2 *The philosophical synthesis.* The philosophy set out in the two pages we have just translated is a synthesis, combining elements of the German and French traditions.

By his own testimony, Kolmogorov drew first and foremost from von Mises. In a footnote, he put the matter this way:

... In laying out the assumptions needed to make probability theory applicable to the world of real events, the author has followed in large measure the model provided by Mr. von Mises ...

The very title of this section of the *Grundbegriffe*, "*Das Verhältnis zur Erfahrungswelt*," echoes the title of the passage in von Mises (1931) that Kolmogorov cites—"*Das Verhältnis der Theorie zur Erfahrungswelt*"—but Kolmogorov does not discuss collectives. As he explained in a letter to Fréchet in 1939, he thought only a finitary version of this concept would reflect experience truthfully, and a finitary version, unlike von Mises' infinitary version, could not be made mathematically rigorous. So for mathematics, one should adopt an axiomatic theory "whose practical value can be deduced directly" from a finitary concept of collectives.

Although collectives are in the background, Kolmogorov starts in a way that echoes Chuprov more than von Mises. He writes, as Chuprov (1910, page 149) did, of a system of conditions (*Komplex von Bedingungen* in German; комплекс условий in Russian). Probability is relative to a system of conditions $\mathfrak{S}$, and yet further conditions must be satisfied in order for events to be assigned a probability under $\mathfrak{S}$. Kolmogorov says nothing more about these conditions, but we may conjecture that he was thinking of the three sources of probabilities mentioned by von Mises: gambling devices, statistical phenomena and physical theory.

Where do von Mises' two axioms—probability as a limit of relative frequency and its invariance under selection of subsequences—appear in Kolmogorov's account? Principle A is obviously a finitary version of von Mises' axiom that identifies probability as the limit of relative frequency. Principle B, on the other hand, is the strong form of Cournot's principle (see Section 2.2.2 above). Is it a finitary version of von Mises' principle of invariance under selection? Evidently. In a collective, von Mises says, we have no way to single out an unusual infinite subsequence. One finitary version of this is that we have no way to single out an unusual single trial. It follows that when we do select a single trial (a single realization of the conditions $\mathfrak{S}$, as Kolmogorov puts it), we should not expect anything unusual. In the special case where the probability is very small, the usual is that the event will not happen.

Of course, Principle B, like Principle A, is only satisfied when there is a collective, that is, under certain conditions. Kolmogorov's insistence on this point is confirmed by the comments we quoted in Section 2.2.2 herein on the importance and nontriviality of the step from "usually" to "in this particular case."

As Borel and Lévy had explained so many times, Principle A can be deduced from Principle B together with Bernoulli's theorem, which is a consequence of the axioms. In the framework that Kolmogorov sets up, however, the deduction requires an additional assumption: we must assume that Principle B applies not only to the probabilities specified for repetitions of conditions $\mathfrak{S}$, but also to the corresponding probabilities (obtaining by assuming independence) for repetitions of $n$-fold repetitions of $\mathfrak{S}$. It is not clear that this additional assumption is appropriate, not only because we might hesitate about independence (see Shiryaev's comments on page 120 of the third Russian edition of the *Grundbegriffe*, published in 1998), but also because the enlargement of our model to $n$-fold repetitions might involve a deterioration in its empirical precision to the extent that we are no longer justified in treating its high-probability predictions as practically certain. Perhaps these considerations justify Kolmogorov's presenting Principle A as an independent principle alongside Principle B rather than as a consequence of it.

Principle A has an independent role in Kolmogorov's story, however, even if we do regard it as a consequence of Principle B together with Bernoulli's theorem, because it comes into play at a point that precedes the adoption of the axioms and hence the derivation of Bernoulli's theorem: it is used to motivate the axioms (cf. Bartlett, 1949). The parallel to the thinking of Lévy is striking. In Lévy's picture, the notion of equally likely cases motivates the axioms, while Cournot's principle links the theory with reality. The most important change Kolmogorov makes in this picture is to replace equally likely cases with frequency; frequency now motivates the axioms, but Cournot's principle remains the most essential link with reality.

In spite of the obvious influence of Borel and Lévy, Kolmogorov cites only von Mises in this section of the *Grundbegriffe*. Philosophical works by Borel and Lévy, along with those by Slutsky and Cantelli, do appear in the *Grundbegriffe*'s bibliography, but their appearance is explained only by a sentence in the preface: "The bibliography gives some recent works that should be of interest from a foundational viewpoint." The emphasis on von Mises may have been motivated in part by political prudence. Whereas Borel and Lévy persisted in speaking of the subjective side of probability, von Mises was an uncompromising frequentist. Whereas Chuprov and Slutsky worked in economics and statistics, von Mises was an applied mathematician, concerned more with aerodynamics than social science, and the relevance of his work on collectives to physics had been established in the Soviet literature by Khinchin (1929; see also Khinchin, 1961, and Siegmund-Schultze, 2004). (For more on the political context, see Blum and Mespoulet, 2003; Lorentz, 2002; Mazliak, 2003; Seneta, 2004.)

5.2.3 *Why was Kolmogorov's philosophy not more influential?* Although Kolmogorov never abandoned his formulation of frequentism, his philosophy has not enjoyed the enduring popularity of his axioms. Section 2 of Chapter I of the *Grundbegriffe* is seldom quoted. Cournot's principle remained popular in Europe during the 1950s (Shafer and Vovk, 2005), but never gained substantial traction in the United States.

The lack of interest in Kolmogorov's philosophy during the past half century can be explained in many ways, but one important factor is the awkwardness of extending it to stochastic processes. The first condition in Kolmogorov's credo is that the system of conditions should be capable of unlimited repetition. When we define a stochastic process in terms of transition probabilities, as in Kolmogorov (1931), this condition may be met, for it may be possible to start a system repeatedly in a given state, but when we focus on probabilities for sets of possible trajectories, we are in a more awkward position. In many applications, there is only one realized trajectory; it is not possible to repeat the experiment to obtain another. Kolmogorov managed to overlook this tension in the *Grundbegriffe*, where he showed how to represent a discrete-time Markov chain in terms of a single probability measure (Chapter I, Section 6), but did not give such representations for continuous stochastic processes. It became more difficult to ignore the tension after Doob and others succeeded in giving such representations.

## 6. CONCLUSION

Seven decades later, the *Grundbegriffe*'s mathematical ideas still set the stage for mathematical probability. Its philosophical ideas, especially Cournot's principle, also remain powerful, even for those who want to go beyond the measure-theoretic framework (Shafer and Vovk, 2001). As we have tried to show in this article, the endurance of these ideas is not due to Kolmogorov's originality. Rather, it is due to the presence of the ideas in the very fabric of the work that came before. The *Grundbegriffe* was a product of its own time.

## ACKNOWLEDGMENTS

Glenn Shafer's research was partially supported by NSF Grant SES-98-19116 to Rutgers University. Vladimir Vovk's research was partially supported by EPSRC Grant GR/R46670/01, BBSRC Grant 111/BIO14428, MRC Grant S505/65 and EU Grant IST-1999-10226 to Royal Holloway, University of London.

We want to thank the many colleagues who have helped us broaden our understanding of the period discussed in this article. Bernard Bru and Oscar Sheynin were particularly helpful. We also benefited from conversation and correspondence with Pierre Crépel, Elyse Gustafson, Sam Kotz, Steffen Lauritzen, Per Martin-Löf, Thierry Martin, Laurent Mazliak, Paul Miranti, Julie Norton, Nell Painter, Goran Peskir, Andrzej Ruszczynski, J. Laurie Snell, Stephen M. Stigler and Jan von Plato.

We are also grateful for help in locating references. Sheynin gave us direct access to his extensive translations. Vladimir V'yugin helped us locate the original text of Kolmogorov's 1929 article, and Aleksandr Shen' gave us a copy of the 1936 Russian translation of the *Grundbegriffe*. Natalie Borisovets, at Rutgers' Dana Library, and Mitchell Brown, at Princeton's Fine Library, have also been exceedingly helpful.